
\input amssym.def
\input amssym.tex
\font\bigbf=cmbx12
%
%
\font\teneusm=eusm10
\font\seveneusm=eusm7
\font\fiveeusm=eusm5
\newfam\eusmfam
\textfont\eusmfam=\teneusm
\scriptfont\eusmfam=\seveneusm
\scriptscriptfont\eusmfam=\fiveeusm

\font\tenmib=cmmib10
\font\sevenmib=cmmib7
\font\fivemib=cmmib5
\newfam\mibfam
\textfont\mibfam=\tenmib
\scriptfont\mibfam=\sevenmib
\scriptscriptfont\mibfam=\fivemib

\font\tenss=cmss10
\font\sevenss=cmss8 scaled 833
\font\fivess=cmr5
\newfam\ssfam
\textfont\ssfam=\tenss
\scriptfont\ssfam=\sevenss
\scriptscriptfont\ssfam=\fivess
\def\ss{\fam\ssfam}
\thinmuskip = 2mu
\medmuskip = 2.5mu plus 1.5mu minus 2.1mu  
\thickmuskip = 4mu plus 6mu
\def\loosegraf#1\par{{%
\baselineskip=13.4pt plus 1pt \lineskiplimit=1pt \lineskip=1.3 pt
#1\par}}

\def\cc{{\Bbb C}}
\def\pp{{\Bbb P}}
\def\eps{\varepsilon}
\def\zz{{\Bbb Z}}
\def\rr{{\Bbb R}}
\def\ss{{\Bbb S}}
\def\bb{{\Bbb B}}

\magnification=\magstep1

\bigskip
\bigskip
\bigskip
\bigskip
\bigskip
\bigskip
\centerline{\bigbf Complex Plateau problem in non K\"ahler manifolds }

\bigskip

\centerline{\rm S. Ivashkovich}
\bigskip
\centerline{December 1997}

\footnote{}{AMS subject classification:
 32 D 15. Key words: meromorphic map, Continuity principle, Hartogs extension
 theorem, spherical shell, complex Plateau problem.}

\bigskip
\bigskip
\centerline{
\vbox{\hsize=4truein
\baselineskip=10pt \sevenrm
In this paper we consider a complex Plateau problem for strong\-ly
pseudoconvex contours in non K\"ahler manifolds. We give the
nessessary and sufficient condition for the existence of solution
in the class of manifolds carrying  pluriclosed metric forms and
propose a conjecture for the general case.
}}

\bigskip
\bigskip
\def\longpoints{\leaders\hbox to 0.5em{\hss.\hss}\hfill \hskip0pt}

\bigskip
\noindent\bf
0. Introduction.\rm
\medskip
 Recall that a complex Plateau problem for a compact real submanifold
$M$ of a complex manifold $X$ consists in finding an analytic chain $A
\subset X\setminus M$ with "boundary" $M$. More preciseliy, one whants to
 find a complex analytic subset $A\subset
X\setminus M$ such that $\partial [A]=[M]$ in the sense of currents. The
nessessary condition on $M$ for complex Plateau problem to have a solution
is the maximal complexity of $M$, i.e. $M$ should be a $CR$-submanifold
of $X$ of ${\sl dim}_{CR} M=p$, where $\dim_{\rr }M=2p+1$. In the case when
$X$ is Stein this is also sufficient, see [H]. Already in the case $X=
\cc\pp^3$ the maximal complexity of the "contour" $M$ is not sufficient
any more, [Db]. We send the interested reader to [Db-H] for an extensive
exposition on the Plateau problem in projective space.

We shall restrict ourselves in this paper by the strongly pseudoconvex
contours $M$, which are already the boundaries of some abstract complex
manifolds. At the same tame we shall look for the solutions of complex
Plateau problem in more general ambient manifolds.

Let $w$ be a strictly positive, smooth $(1,1)$-form on the complex
manifold $X$.

\smallskip\noindent\bf
Definition 0.1. \it We shall say that $w$ is pluriclosed if $dd^cw=0$.

\smallskip\rm
A Hermitian metric canonically associated to a such $w$ we shall often
also call \it pluriclosed\rm .

 \smallskip\noindent\bf
 Definition 0.2. \it We say that a complex space $X$ is disk-convex in
 dimension $k$ if for any compact $K\subset\subset X$ there is another
 compact $\hat K$ such that for every meromorphic mapping $f:\bar\Delta^k
 \to X$ with $f(\partial \Delta^k)\subset K$ one has
 $f(\bar\Delta^k)\subset \hat K$.
 \smallskip\rm
 All compact spaces are of course disk-convex. More generally all $k$-
 convex spaces are disk-convex in dimension $k$.

 \smallskip\noindent\bf
 Definition 0.3. \it A maximally complex $CR$-manifold $M$ we call strongly
 pseudoconvex if it can be realized as a strongly pseudoconvex hypersurface
 in some complex manifold.

 \smallskip\rm
 If we say that a $CR$-manifold $M$ is contained in a complex manifold $X$
 we mean that some $CR$-imbedding $M\to X$ is given.

\smallskip\noindent\bf
Theorem. \it Let $M$ be a strongly pseudoconvex, maximally complex,
compact $CR$-manifold in a disk-convex complex manifold $X$ carrying
a pluriclosed Hermitian metric form. Suppose that $M$ bounds an abstract
smooth Stein domain.

(a) If $\dim M\ge 5$ then the complex Plateau problem for $M\subset X$ has a
solution.

(b) If $\dim M=3$ then the complex Plateau problem for $M\subset X$ has a
solution iff $M$ is

homologous to zero in $X$.
\smallskip\noindent\bf
Remarks. \rm 1. Let $H^2:=\cc^2\setminus \{ 0\} /(z\sim 2z)$ be a Hopf surface.
Take $M$ to be an image of a standard unit sphere from $\cc^2$ under the
natural projection $\pi :\cc^2\setminus \{ 0\} \to H^2$. $M$ is not homologous
to zero in $H^2$ (i.e. it is a spherical shell in $H^2$!), so a complex Plato
problem has no solution for this $M$. Note that $w={i\over 2}{dz_1\wedge
d\bar z_1 + dz_2\wedge d\bar z_2\over \Vert z\Vert ^2}$ is a pluriclosed
hermitian metric form on $H^2$.

In the case (a), i.e. when $\dim M\ge 5$ the spherical shells in $X$ are
not an obstructions for finding a film with boundary $M$ because we have
"enough concavity".

\noindent 2. Consider a Hopf three-fold  $H^3:=\cc^3\setminus \{ 0\} /(z\sim
2z)$.
In this case take a sphere $\ss^3$ in a hyperplane $\{ z_1=0\} $. Its image
$M$ under the natural projection will be homologous to zero but will not
bound any analytic set in $H^3$. The reason here is that $H^3$ doesn't
admits a pluriclosed Hermitian metric.

\noindent
3. If one doesn't requires strict pseudoconvexity of a "contour" $M$ then
counterexamples are known already in $\cc\pp^3$, see [Db].

\noindent 4. The condition on $M$ to bound an abstract smooth Stein
domain is really restrictive in dimension 3, while for higher dimensions
one has the Rossi theorem guaranteeing the existence (but in general not
smooth) abstract Stein domain with boundary $M$, see [Rs].

\noindent 5. When $X$ is K\"ahler then \it any  str.ps.convex $M\subset X$,
\rm which bounds an abstract Stein domain, is homologous to zero in $X$. This
follows from the Hartogs-type extension theorem proved [Iv-1].

The proof of the Theorem consists in meromorhpic extension of a
$CR$-imbedding of $M$ into $X$ onto the
Stein domain bounded by this $M$. We do it along the levels of appropriate
plurisubharmonic Morse exhaustion function, see \S\S 2 and 3.
\smallskip\noindent\sl
\smallskip\rm
\bigskip\noindent\bf
Table of content\medskip \sl
\noindent
\line{0. Introduction. \longpoints pp. 1-2}

\smallskip\noindent
\line{1. Continuity principle and spherical shells. \longpoints pp. 2 - 5}

\smallskip\noindent
\line{2. Three-dimensional contours. \longpoints pp. 5 -8}

\smallskip\noindent
\line{3. Higher-dimensional contours. \longpoints pp. 8 - 9}

\line{4. Open questions. \longpoints pp. 9 - 10}

\smallskip\noindent
\smallskip\noindent
\line{References. \longpoints pp. 10 - 12}

 \medskip\noindent\bf
 1. Continuity principle and spherical shells.
 \smallskip\rm
 To realise this approach we shall need some resuts on extension of
 meromorphic mappings into a general complex manifolds (spaces). For the
 convenience of the reader we collect them in this paragraph.

 Put $A^k(r,1)=\{ z\in \cc^k: r<\Vert z\Vert <1\} $ and $A^k_s(r,1):=\{ s\}
 \times A^k(r,1)$ for $s\in \Delta^n$. Let $f:\Delta^n\times A^k(r,1)\to X$
 be a holomorphic mapping into a normal complex space $X$. Denote by $S$
 the set of points $s\in \Delta^n$ such that the restriction
 $f_s:=f\mid_{A^k_s(r,1)}$ extends meromorphically onto the polydisk
 $\Delta^k_s:=\{ s\}\times \Delta^k$.

 \smallskip\noindent\bf
 Theorem 1.1. \sl (Continuity principle). \it
  Let $f:\Delta^n\times A^k(r,1)\to X$
 be a holomorphic mapping into a normal, disk-convex in dimension $k$
 complex space $X$.
 Suppose that there is a constant $C_0<\infty$ and a compact $K\subset X$
 such that for $s$ in some
 subset $S\subset \Delta^n$, which is thick at origin:

 (a) the restriction  $f_s:=f\mid_{A^k_s(r,1)}$ is well defined and extends
 meromorphically
 onto the polydisk  $\Delta^k_s:=\{ s\}\times \Delta^k$,  and
 ${\sl Vol}(\Gamma_{f_s})\le C_0$ for all $s\in S$;

 (b) $f(\Delta^n\times A^k(r,1))\subset K$ and $f_s(\Delta^k)\subset K$
 for all $s\in S$.

\noindent
Then:

1. If $n=1$ then there is a neighborhood $U\ni 0$ in $\Delta $ such that
 $f$ extends meromorphically onto $U\times \Delta^k$.

2. If $n\ge 2$ and $X$ has bounded cycle geometry in dimension $k$, then
again there is a neighborhood $U\ni 0$ in $\Delta^n$ and a meromorphic
extension of $f$ onto $U\times \Delta^k$.

 \smallskip\rm Here $\Gamma_{f_s}$ denotes the graph of $f_s$ and the volumes
 are taken with respect to some Hermitian  metric $h$ on $X$
 and a standart Euclidean metric on $\cc^k$. The condition of finitness clearly
 doesn't depends on the particular choice of the metrics. As usually saying
 that a set $S\subset \Delta^n$ is thick at the point $z_0$  we mean that
 for any
 neighborhood $U\ni z_0$ $S\cap U$ is not contained in a proper analytic
 subset of $U$.

 Denote by ${\cal B}_k(X)$ the Barlet space of compact analytic cycles
 of dimension $k$ in $X$. This is an analytic space, which has not more
 then countable number of components.

 \smallskip\noindent\bf
 Definition 1.1. \it We say that a complex space $X$ has bounded cycle
 geometry in dimension $k$ if all connected components of the Barlet space
 ${\cal B}_k(X)$ are compact.

 \rm In other words $X$ has bounde cycle geometry in dimension $k$ if all
 \it irreducible \rm components of ${\cal B}_k(X)$ are compact and all
 connected components of ${\cal B}_k(X)$ are just the finite unions of
 irreducible ones.
 Note again
 that the property to have bounded cycle geometry doesn't depend on the
 choice of Hermitian metric.

 For the proof of this theorem we refer to [Iv-2]. There also an example
 is given, showing that in the case $n\ge 2$ some condition on $X$
 is needed on the contrary to the case $n=1$.

 There are some cases, occuring in the applications, when the condition
 of bounded cycle geometry can be dropped. For example one  has the
 following

   \smallskip\noindent\bf
   Proposition 1.2. \it Let $f:\Delta^n\times A^k(r,1)\to X$ be a holomorphic
   map into a normal, disk-convex in dimension $k$ complex space $X$.
   Suppose that:

   (1) for every $s\in \Delta^n$ outside of thin set, the restriction
   $f_s$ extends meromorphically onto $\Delta_s^k$;

   (2) there is a compact $K\subset\subset X$ such that $f_s(\Delta_s^k)
   \subset K$ for all $s$ and $f(\Delta^n\times A^k(r,1))\subset K$;

   (3) the volumes of the graphs $\Gamma_{f_s}$ are uniformly bounded in
   $\Delta^n$, i.e. there exists $C_0<\infty $ s.t.
   ${\sl vol}(\Gamma_{f_s})\le C_0$ for all $s$.

   \noindent
   Then $f$ meromorphically extend onto $\Delta^{n+k}$.
   \smallskip\noindent\rm
   For the proof see [Iv-2], Corollary 1.

   It is worth, probably to point out one case when the boundedness
   of cycle geometry is satisfied automatically - when $k={\sl dim}X - 1$.
   Really the cycle space of divisors is allwayse compact (provided $X$
   is compact).

  Denote by

 $$
 H^k_n(r):=\{ (z',z^{''})\in \Delta^{n+k}:1-r<\Vert z{''}\Vert <1 \hbox{or}
 \Vert z'\Vert <r\} =
 $$
 $$
 = \Delta^n\times A^k(1-r,1)\cup \Delta^n(r)\times \Delta^k\eqno(1.1)
 $$
 the $k$-concave Hartogs figure in $\cc^{n+k}$.
 \smallskip\noindent\bf
 Definition 1.2. \it We say that the meromorphic mappings into the space $X$
  have the
 Hartogs-type extension property in bidimension $(n,k)$ if any meromorphic
 map $f:H^k_n(r)\to X$ extends meromorphically onto $\Delta^{n+k}$.

 \smallskip\noindent\bf
 Definition 1.3. \it Let us call a Hermitian form $w$ on $X$ plurinegative if
 $dd^cw\le 0$.
 \smallskip\rm
 The class of normal complex spaces admiting plurinegative Hermitian metrik
 form we shall denote by ${\cal P}_{-}$.
 \smallskip\noindent\bf
 Theorem 1.3. \it Let $f:H^1_n(r)\to X$ be a meromorphic map into a disk-convex
  complex space $X$ which admits a plurinegative Hermitian metrik form. Then:

  (1) $f$ extends to a meromorphic map  $\hat f:\Delta^{n+1}\setminus A
  \to X$, where $A$ is closed $(n-1)$-polar subset of $\Delta^{n+1}$.

  (2) If moreover, $w$ is pluriclosed then $A$ is a analytic subvariety of
  $\Delta^{n+1}$ of pure codimension two (may be empty). If $A\not= \emptyset
   $ then for every sphere $\ss^3$ embedded into $\Delta^{n+1}
  \setminus A$ in such a way that $[\ss^3]\not=0$ in $H_3(\Delta^
  {n+1}\setminus A,\zz)$, its image $f(\ss^{3})$ also is not homologous to zero
  in $X$.

\smallskip\noindent\bf
Remarks. \rm 1. One can estimate the number of irreducible components of
the singularity set $A$ in this theorem meating a  compact subset $P\subset
\subset \Delta^{n+1}$. Namely, let a compact $K\subset X$, which contains
${\sl cl}[f(P\setminus S)]$, is chosen to be a finite  subcomplex of $CW$-
complex $X$. Choose a point $z'\in \Delta^{n-1}$ such that $A$ intersects
$\Delta^2_z:=\{ z\} \times \Delta^2$ by discrete set $A_{z'}$. Let $A_{z'}
\cap \partial P=\emptyset $. Then

$$
\mid A_{z'}\cap P\mid \le \vert \int_{\partial (P\cap \Delta^2_{z'})}
d^cw\vert \cdot [{\sl inf}\{ \vert \int_{\gamma }d^cw\vert :\gamma \in H_3(K,\zz ),
\int_{\gamma }d^cw\not=0\} ]^{-1}.\eqno(1.2)
$$

In other words the number of branches of singular set (and moreover, their
existence)
is bounded by the differential geometry of $X$. Remark that the subset $\{
\vert \int_{\gamma }d^cw
\vert : \gamma \in H_3(K,\zz ), \int_{\gamma }d^cw\not=0\}\subset \rr $ is
separated from zero, see (2.2.14) in [Iv-2].

\noindent
2. Let us call a spherical shell of dimension $k+1$ in complex
space
$X$  an image $\Sigma $ of the standard sphere $\ss^{2k+1}\subset
\cc^{k+1}$
under the meromorphic map of some neighborhood of $\ss^{2k+1}$ into $X$,
such
that $\Sigma $ is not homologous to zero in $X$. This notion is close to
the notion of the global spherical shell, introduced by Kato, see [Ka-3].
Thus we obtain the following
\smallskip\noindent\bf
Corollary 1.4. \it Let $X$ be a disk-convex complex space which possess a
pluriclosed Hermitian metric form. Then the following is equivalent:

 (a) $X$ possesses a meromorphic extension property in bidimension $(n,1)$
 for all $n\ge 1$,

 and thus in all bidimensions $(n,k)$.

 (b) $X$ contains no two-dimensional spherical shells.

\smallskip\noindent\rm
3. A wide class of complex manifolds without two-dimensional
spherical shells is  for example a class of such manifolds $X$ for which the
Hurewicz homomorphism  $\pi_3(X)\to H_3(X,Z)$ vanishes.

\smallskip
In the proof of the results listed in this paragraph, as well as for the
solution of the complex Plateau problem in this paper, we use the following
lemma. Consider a meromorphic mapping $f:\Delta^p\times \Delta^q(a)\to X$
into a complex space $X$. Here $\Delta^q(a)$ is a polydisk in $\cc^q$ of
radii $a$, $\Delta^q(1)=\Delta^q$. Let $S$ be some closed subset
of $\Delta^p$ and $s_0\in S$. Suppose that for each $s\in S$ the restriction
$f_s:=f\mid_{\{ s \} \times \Delta^q(a)}$ is well defined and meromorphically
extends onto a $q$-disk $\Delta^q(b)$, $b>a$. Denote by $\nu_j=\nu_j(K)$ the
minima
of volumes of $j$-dimensional compact analytic subsets contained in some
compact $K\subset X$, see Lemma 2.3.1 from [Iv-3]. Fix some $a<c<b$. Put
$$
\nu = min\{ vol(A_{q-j}\cdot \nu_j:j=1,...,q\} , \eqno(1.3)
$$
\noindent
where $A_{q_j}$ are running over all $(q-j)$-dimensional analytic subsets of
$\Delta^q(b)$, intersecting $\Delta^q(c)$.

\smallskip\noindent\bf
Lemma 1.5. \it Suppose that there exists a neighbourhood $U\ni s_0$ in
$\Delta^p$ such that for all $s_1,s_2\in S\cap U$
$$
\vert vol(\Gamma_{f_{s_1}}) - vol(\Gamma_{f_{s_2}})\vert <\nu /2, \eqno(1.4)
$$
\noindent and that $s_0$ is a locally regular point of $S$. Then there
exists a neighbourhood $V_c\ni s_0$ in $\Delta^p$, such that $f$
meromorphically exteds onto $V_c\times \Delta^q(c)$.

\smallskip\rm For the proof see [Iv-3], Lemma 2.4.1.

\magnification=\magstep1

\bigskip\noindent\bf
2. Three-dimensional contours.
\smallskip\rm

We shall prove in this paragraph the part (b) of our {\sl Theorem} from
{\sl Introduction}.

 We suppose that $M$ bounds an abstract smooth Stein domain, i.e. there is a
complex manifold $D $ with boundary $M$ such that $D \setminus
M$ is Stein, and that a $CR$-imbedding $f:M\to X$ is given. All that we  need
to prove is that $f$ extends meromorphically onto $D $. Clearly we
can suppose that $f$ is already holomorphically extended to some neighborhood
of $M$ in $D $.
\smallskip\noindent\bf
Proposition 2.1. \it Let $(D, M)$ be as above and suppose additionally
that $\dim M=3$.

(a) Then any $CR$-map $f:M\to X$, where $X$ is a disk-convex complex space
admitting

a pluriclosed Hermitian metric form, extends meromorphically onto
$D \setminus S$. Here $S$

is a finite subset of $D$.

(b) If $f(M)$ is homologous to zero in $X$, or if $X$ doesn't contain
spherical shells, then

$S$ is empty.
\smallskip\noindent\sl
Proof. \rm Let $\rho :D \to [0,1]$ be a strictly plurisubharmonic (and thus
Morse) exhausting function. Denote by $D^{+}_{\eps }=\{ z\in D: \rho (z)>
\eps \} $.
Let ${\cal E}$ be the set of such $\eps $ that $f$ can be meromorphically
extended onto $D^{*}_{\eps }\setminus S_{\eps }$, where $S_{\eps }$ is a
discrete set. ${\cal E}$ is obviously closed and nonempty.
All we need to prove is that ${\cal E}$ is open.

Let $\eps_0 =\inf \{ \eps \in {\cal E}\} $. If $\eps_0$ is a regular value
of $\rho $ then the needed result immediately follows from part (2) of
{\sl Theorem 1.3}.

Consider the case of not regular value $\eps_0$ of $\rho $. Denote by
$M_{\eps_0}=\{ z:\rho (z)=\eps_0\} $-the critical level set. Fix a critical
point
$z_0\in M_{\eps_0}$. All we need to prove is that for any neighborhood
$W$ of $z_0$ the envelope of holomorphy of $W\cap D^{+}_{\eps }$ contains
some neighborhood of $z_0$. For convenience we can suppose that $z_0=0$ and
$\eps_0=0$. Write
$$
\rho (z) = Q(z) + <z,z> + \bar Q(z) + O(\Vert z\Vert^3), \eqno(2.1)
$$
\noindent where $Q(z)$ is a holomorphic polynomial, $<z,z>$ - Hermitian form -
 Levi form of $\rho $. By linear coordinate change we transform $< , >$ to the
sum of squares of absolute values. Then by unitary coordinate change we transform
$Q$ to the some of squares with real nonnegative coefficients. Now (2.1)
has a form
$$
\rho (z)= \sum_{j=1}^pa_jz_j^2 + \sum_{j=1}^pa_j\bar z_j^2 + \sum_{j=1}^n
\vert z_j\vert^2 + O(\Vert z\Vert^3). \eqno(2.2)
$$
\noindent In coordinates $z_j=x_j+iy_j$ we rewrite (2.2) as follows
$$
\rho (z) = 2\sum_{j=1}^pa_j(x_j^2-y_j^2) + \sum_{j=1}^n(x_j^2+y_j^2) +
O(\Vert z\Vert^3) =
$$
$$
= \sum_{j=1}^p[(1+2a_j)x_j^2 + (1-2a_j)y_j^2] + \sum_{j=p+1}^n(x_j^2+y_j^2)
+ O(\Vert z\Vert^3). \eqno(2.3)
$$
\noindent Renumerate the coordinates in such a way that $a_j\ge {1\over 2}$
for $j=1,...,q$ and $a_j<1/2$ for $j=q+1,...,p$. Then
$$
\rho (z) \ge \sum_{j=1}^q[(2a_j+1)x_j^2-(2a_j-1)y_j^2] + \delta\cdot \sum_
{j=q+1}^p\vert z_j\vert^2 + O(\Vert z\Vert^3) \ge
$$
$$
\ge \sum_{j=}^q[(2a_j-\delta_1+1)x_j^2-(2a_j+\delta_1-1)y_j^2] + \delta\cdot
\sum_{j=q+1}^p\vert z_j\vert^2  := \rho_1(z), \eqno(2.4)
$$
\noindent for some $\delta >0$ and $\delta_1$ can be chosen arbitrarily small
for small $\Vert z\Vert $. While obviously $D^+:=\{ z\in \bb^n: \rho_1(z)>0\}
\subset D^+_{\eps_0}$, all we need is to prove the following
\smallskip\noindent\sl
Lemma 2.2. \it The envelope of holomorphy of $D^+$ contains the origin.
\smallskip\noindent\sl
Proof. \rm Consider two cases.

\noindent\sl
Case 1: $q\le n-1$. \rm In this case $D^+$ contains the following Hartogs
figure:
$$
H := \{ z\in \bb^n: \sum_{j=1}^q[(2a_j-\delta_1+1)x_j^2 - (2a_j+\delta_1-1)
y_j^2]>0, \delta \cdot \sum_{j=q+1}^n\vert z_j\vert^2<1
$$
or
$$
\sum_{j=1}^q(2a_j-\delta_1+1)x_j^2 - (2a_j+\delta_1-1)y_j^2]>-\eps , \delta
\cdot \sum_{j=q+1}^n\vert z_j\vert^2>\eps \} .
$$
\noindent
The envelope of holomorphy of $H$ obviously contains the origin.
\smallskip\noindent\sl
Case 2: $q=n$. \rm In this case
$$
D^+ = \{ z\in \bb^n: \sum_{j=1}^n[(2a_j-\delta_1+1)x_j^2 - (2a_j+\delta_1-1)
y_j^2]>0\} .\eqno(2.5)
$$
\noindent Put $b_j=2a_j-\delta_1+1$, $c_j=2a_j+\delta_1-1$, $j=1,...,n$. For
small $\delta-1$, $b_j>c_j$. Write (2.5.5) in the form
$$
D^+ = \{ z\in \bb^n: \sum_{j=1}^nb_jx_j^2 >  \sum_{j=1}^nc_jy_j^2\} .
\eqno(2.6)
$$
\noindent In the new coordinates $z_j\to \sqrt{b_j}z_j$ (2.6) take a form
$$
D^+ = \{ z\in \bb^n:  \sum_{j=1}^nx_j^2 >  \sum_{j=1}^n\delta_jy_j^2\} ,
\eqno(2.7)
$$
\noindent where $\delta_j = {c_j\over b_j}<1$, $j=1,...,n$. Put $\delta_0:=
\max \{ \delta_1,...,\delta_n\} <1$. Then
$$
D^+\supset D^+_1=\{ z\in \bb^n: \Vert x\Vert^2>\delta_0\cdot \Vert y\Vert^2\}
.\eqno(2.8)
$$
\noindent The set $D_1^+$ contains clearly the following complete "tube torus"
$$
T = \{ x+iy\i  \cc^n: \Vert x\Vert =1, \Vert y\Vert \le 1/\delta_0\} ,
\eqno(2.9)
$$
\noindent
where $1/\delta_0 := \eta >1$. We shall prove that already the envelope of
holomorphy of $T$ contains the origin. For this consider the following
continuous family of complex hypersurfaces
$$
C_t = \{ z\in \cc^n: z_1^2+...+z_n^2=t\} \eqno(2.10)
$$
\noindent or
$$
C_t = \{ x+iy\in \cc^n:  \Vert x\Vert^2-\Vert y\Vert^2=t, (x,y)=0 \} ,
\eqno(2.11)
$$
\noindent where $(x,y)=x_1y_1+...+x_ny_n$. Consider the intersections of $C_t$
with a ball of radii $1+\eta^2$:
$$
\tilde C_t = \{ x+iy \in \bb^n_{1+\eta^2 }: \Vert x\Vert -\Vert y\Vert =t,
(x,y)=0\} .\eqno(2.12)
$$
\noindent This is a continuous family of irreducible analytic hypersurfaces in
$\bb^n_{1+\eta }$ such that
$$
\tilde C_{1+\eta^2 } = \{ x+iy\in \bb^n_{1+\eta^2 }: \vert x\Vert^2-\Vert y
\Vert^2=1+\eta^2, (x,y)=0\} =
$$
$$
= \{ x+iy\in \bb^n_{1+\eta^2}: \Vert x\Vert^2+\Vert y\Vert^2=1+\eta^2=\Vert x
\Vert^2-\Vert \Vert^2, (x,y)=0\} =
$$
$$
\{ x+iy\in \bb^n_{1+\eta^2}: \Vert x\Vert^2=1+\eta^2, y=0\} \subset T,
$$
\noindent but $\tilde M_0\ni 0$. By continuity principle the envelope of
holomorphy of $T$ contains the origin.
\smallskip
\hfill{q.e.d}
\smallskip\noindent\sl
End of the proof. \rm

So, as in the case of regular value, we can extend our map $f$ meromorphically
to the neighborhood of the critical level $M_{\eps_0}$ minus discrete set.
 As a result we obtain the extension $\hat f$ of our map onto $\bar D\setminus
S$ where $S$ is a finite subset of $\bar D$ not intersecting $M=\partial D$.
 If we put $T:= f^*w$ then $dd^c\tilde T$ is nonpositive measure supported on
$S$, see Lemma 2.6.1 from [Iv-2]. We have
$$
\int_Sdd^c\tilde T = \int_Ddd^c\tilde T = \int_{\partial D}d^cT = \int_{f(
\partial D)}d^cw = \int_Md^cw=0,
$$ \noindent if $M$ is homologous to zero in $X$, or if $X$ doesn't contain
spherical shells.
\smallskip
\hfill{Part (b) of the Theorem is proved.}

\bigskip\noindent\bf
3. Higher-dimensional contours.
\smallskip\rm
 We shall prove now the part (a) of our {\sl Theorem} from {\sl Introduction}.
In the proof we shall need the following statement. Let $S$ be a closed
subset in the product of two balls $W=B^{n-1}\times B^2$ such that
$S\cap (\overline{B^{n-1}}\times \partial B^2)=\emptyset $. For $z^{'}\in
B^{n-1}$ denote $B^2_{z^{'}}:=\{ z^{'}\} \times B^2$. Let a meromorphic
mapping mapping $f:W\setminus S\to X$ is given, where $X$ carries a
pluriclosed metric form $w$. Suppose that $f$ is holomorphic in the
neighborhood of $\overline{B^{n-1}}\times \partial B^2$. Denote by
$T$ the preimage of $w$ by $f$ and suppose that $T$ admits a trivial
extension $\tilde T$ onto $W$.

\smallskip\noindent\bf
Lemma 3.1. \it Suppose that for all $z^{'}\in B^{n-1}$ $f(B^2_{z^{'}})$ is
homologous to zero in $X$. Then:

(i) $dd^c\tilde T=0$ in the sence of currents.

(ii) There is a $(1,0)$-current $\gamma $ on $W$, smooth in the neighbourhood
of $\overline B^{n-1}\times \partial B^2$, such that $\tilde T=i(\partial
\bar\gamma - \bar\partial\gamma )$.

\smallskip\noindent\rm
For the proof see [Iv-2], Lemma 2.6.2.
\smallskip
Remark that  only at the end of the proof in \S 2 we used the fact that the
dimension of $D$ is two. So the following proposition clearly enables us to
finish the proof also of  part (a) of the {\sl Theorem} i.e for $\dim D\ge 3$.
\smallskip\noindent\bf
Proposition 3.2. \it Every holomorphic map $f$ from $H_1^n(r)$ into a
disk-convex
complex space $X$, which admits a pluriclosed Hermitian metric, extends
meromorphically onto $\Delta^n$
provided $n\ge 2$.
\smallskip\noindent\sl
Proof. \rm It will be convenient for us simultaneously with the proof of the
main statement of the {\sl Proposition} to prove also the following weaker
statement. Denote by $A^n(a,b):=\Delta^n(b)\setminus \Delta^n(a)$, for $0\le
a<b$.

\smallskip
\it Every holomorphic map $f:A^n({1\over 2},1)\to X$, where $X$ from our
{\sl Proposition },

extends meromorphically onto $\Delta^n$, provided $n\ge 2$ and $f(\partial
\Delta^n_{3/4})$ is

homologous to zero in $X$.
\smallskip\rm
We shall prove both statements by  induction on $n$. For $n=2$ the
second statement follows directly from {\sl Theorem 1.1}. So it is sufficient
to prove that for
any $n\ge 2$ from the second statement follows the statement of {\sl
Proposition } for this $n$.

 So let a holomorphic mapping $f:H_1^n(r)\to X$ is given. For every $z\in
\Delta $ restriction $f_z$ of $f$ onto $\Delta^n_z:=\{ z\} \times \Delta^n$ is
holomorphic on $A^n(r,1)$. So, by the assumption $f_z$ meromorphically
extends onto $\Delta^n$, because $f(\partial\Delta^n_z)\sim f(\partial\Delta^
n_0)\sim 0$ in $X$! Lemma 1.5
immediately gives us (after shrinking $\Delta^{n+1}$ and taking different bends
of $z_2,...,z_{n+1}$) the meromorphic extension of $f$ onto $\Delta^{n+1}
\setminus S$. Where $S$ is zerodimensional pluripolar compact in $\Delta^{n+1}
$.

Because $I(f)$ is an analytic set of positive
dimension outside of zerodimensional set, $\overline{I(f)}$ is analytic
in $\Delta^{n+1}\setminus H_1^n$, and thus empty. So the fundamental set
of $f$ is discrete in $\Delta^{n+1}\setminus S$.

Put $T:=f^*w$, where $w$ is pluriclosed metric form on $X$. $T$ has locally
summable coefficients in $\Delta^{n+1}$ and its trivial extension $\tilde T$
is plurinegative with $dd^c\tilde T$ supported on $S$. Observe that
$\tilde T=T$ is pluriclosed outside of $S$.

{\sl Lemma 3.1} tells us that $dd^c\tilde T=0$ and moreover there is a
(1,0)-current $\gamma $ in any given ball $W\subset \Delta^{n+1}$, $\partial
W\cap S=\emptyset $, smooth on $W\setminus S$, such that $\tilde T=i(\partial
\bar\gamma -\bar\partial \gamma )$ . Remark that the conditions of {\sl Lemma
3.1} are satisfied, because $S$ is zerodimensional and $n+1\ge 3$. All that
remained is to repeat the arguments from the proof of {\sl Lemma 2.6.3} to
estimate the volume of the graph of $f$ in the neighborhood of $S$. Namely

$$
{\sl Vol}(\Gamma_{f\mid_{W\setminus S}}) = \int_{W\setminus S}(T + dd^c\Vert
z\Vert^2)^{n+1} = \sum_{j=0}^{n+1}C_{n+1}^j\int_{W\setminus S} T^j\wedge
(dd^c\Vert z'\Vert^2)^{n+1-j}\le
$$

$$
\le C\cdot \int_{W\setminus S} T^j\wedge (dd^c\Vert
z'\Vert^2)^{n+1-j} = C\cdot \lim_{\eps \searrow 0}\int_{W\setminus S}\tilde
T^j_{\eps }\wedge (dd^c\Vert z'\Vert^2)^{n+1-j}\le
$$

$$
\le C\cdot  \lim_{\eps \searrow 0}\int_W\tilde T^j_{\eps }\wedge (dd^c\Vert z'
\Vert^2)^{n+1-j} = C\cdot \lim_{\eps \searrow 0}\int_W(\partial \bar\gamma^
{1,0}_{\eps }+\bar\partial \gamma^{1,0}_{\eps })^j\wedge
(dd^c\Vert z'\Vert^2)^{n+1-j} =
$$

$$
= C\cdot \lim_{\eps \searrow 0}\int_W(d(\bar\gamma^{1,0}_{\eps }+\gamma^{1,0}
_{\eps }))^j\wedge (dd^c\Vert z'\Vert^2)^{n+1-j} =
$$

$$
 = C\cdot \lim_{\eps
\searrow 0}
\int_{\partial W}(\bar\gamma_{\eps }^{1,0}+\gamma_{\eps }^{1,0})\wedge
d(\bar\gamma^{1,0}_{\eps }+\gamma^{1,0}_{\eps })^{j-1}\wedge (dd^c\Vert z'
\Vert^2)^{n+1-j} =
$$

$$
 = C\cdot \int_{\partial W}(\bar\gamma^{1,0}+\gamma^{1,0})\wedge d(\bar\gamma^
{1,0}+\gamma^{1,0})^{j-1}\wedge (dd^c\Vert z'\Vert^2)^{n+1-j} <  \infty
$$
\smallskip
 From Bishop theorem we get an extension of the graph of $f$ onto $\Delta
^{n+1}$.
\smallskip
\hfill{q.e.d.}
\medskip\noindent\bf
4. Open questions.
\smallskip\rm
 \rm Let us propose some open questions arising naturally from the
 exposition.
  Our point of depart will be the following observation due to Gauduchon:

 \it every compact complex manifold of dimension $k+1$ carries a Hermitian

 metric form $w$ with $dd^cw^k=0$.

 \rm Really, the
 condition to carry $dd^c$-closed strictly positive $(k,k)$-form for a
 compact complex
 manifold is alternative to that of carrying a   bidimension $(k+1,
 k+1)$-current $T$ with $dd^cT\ge 0$ but $\not \equiv 0$. This in the case of
 ${\sl dim}X=k+1$ is a nonconstant plurisubharmonic function, which on
 compact $X$ doesn't exist. In fact in [Gd] a stronger statement was
 proved, but we shall not need it here.

  Let us introduce the class ${\cal G}_k$ of normal complex spaces,
 carrying
 a nondegenerate positive $dd^c$-closed strictly positive $(k,k)$-forms.
 Note that the sequence
 $\{ {\cal G}_k\}$ is rather exaustive: ${\cal G}_k$ contains all compact
 complex
 manifolds of dimension $k+1$.

 Note also that compact spaces from ${\cal G}_k$
 have bounded cycle geometry in dimension $k$, see 1.4 from [Iv-2]. We
 conjecture that
 meromorphic mappings into the spaces of class ${\cal G}_k$ are ''almost
 Hartogs-extendable'' in bidimension $(n,k)$ for all $n\ge 1$:

 \smallskip\noindent\sl
 Conjecture 1. \it Every meromorphic map $f:H^k_n(r)\to X$, where $X\in
 {\cal G}_k$ and is disk-convex in dimension $k$, extends to a meromorphic
 map from $\Delta^{n+k}\setminus A$
 to $X$,
 where $A$ is an analytic subvariety of $\Delta^{n+k}$ (may be empty) of pure
 codimension
 $k+1$. Moreover, if $A\not=\emptyset $, then for every sphere $\ss^{2k+1}$
 embedded into $\Delta^{n+k}\setminus A$ in such a way that $[\ss^{2k+1}]
 \not=0$ in $H_{2k+1}(\Delta^{n+k}\setminus A,\zz)$, $f(\ss^{2k+1})$
 also is not homologous to zero in $X$.
 \smallskip\rm
 In [Iv-2] this conjecture is proved for  the case $k=1$.

 \smallskip\noindent\sl
 Conjecture 2. \it Let $M$ be a strictly pseudoconvex, compact contour of
 real dimension $2k+1$ in the $k$-disk-convex complex manifold
 $X\in {\cal G}_k$, $k\ge 2$. Then the Plateau problem for $M$ has
 solution iff $M$ is homologous to zero in $X$.

 \smallskip\rm
 This conjecture would follow from the first one, but is probably easier.
 \smallskip\noindent\sl
 Conjecture 3. \it Let $M$ be a strictly pseudoconvex, compact, three
 dimensional  contour in  compact manifold $X$. Prove that
 $M$ bounds an abstract Stein domain.

 \smallskip\noindent
\magnification=\magstep1
\spaceskip=4pt plus3.5pt minus 1.5pt
\spaceskip=5pt plus4pt minus 2pt
\font\csc=cmcsc10
\font\tenmsb=msbm10
\def\rr{\hbox{\tenmsb R}}
\def\cc{\hbox{\tenmsb C}}
\newdimen\length
\newdimen\lleftskip
\lleftskip=2.5\parindent
\length=\hsize \advance\length-\lleftskip
\def\entry#1#2#3#4\par{\parshape=2  0pt  \hsize%
\lleftskip \length%
\noindent\hbox to \lleftskip%
{\bf[#1]\hfill}{\csc{#2 }}{\sl{#3}}#4%
\medskip
}
\ifx \twelvebf\undefined \font\twelvebf=cmbx12\fi

\bigskip\bigskip
\bigskip\bigskip
\centerline{\twelvebf References.}
\bigskip



\entry{Db}{Dolbeault P.:}{Geometric measure theory and the calculus of
variations.} Proc.\ Symp.\ Pure\ Math. {\bf44}, 171-205, (1986).

\entry{Dl-H}{Dolbeault P., Henkin G.:}{Surfaces de Riemann de bord donne
dans $\cc\pp^n$.} In Contribution to complex analysis and analytic
geometry, ed. by H. Skoda et al. \ Vieweg, Aspects Math. E 26, 163-187
 (1994).

\entry{Ga}{Gauduchon P.:}{Les metriques standard d'une surface a premier
nombre de Betti  pair.} Asterisque.\ Soc.\ Math.\ France.  {\bf126}, 129-135,
(1985).

\entry{H}{Harvey R.:}{Holomorphic chains and their boundaries.} Proc.\ Symp.\
Pure\ Math. {\bf30}, Part I, 307-382 (1977).

\entry{H-L}{Harvey R., Lawson H.:}{An intrinsic characterisation of K\"ahler
manifolds.} Invent.\ math. {\bf74}, 169-198, (1983).

\entry{Iv-1}{Ivashkovich S.:}{The Hartogs-type extension theorem for the
meromorphic maps into compact K\"ahler manifolds.} Invent.\ math. {\bf109}
 , 47-54, (1992).

\entry{Iv-2}{Ivashkovich S.:}{Continuity principle and extension properties
of meromorphic mappings with values in non K\"ahler manifolds. }
MSRI\ Preprint\ No.\ 1997-033.

\entry{Iv-3}{Ivashkovich S.:}{One example in concern with extension and
separate analyticity properties of meromorphic mappings.} Preprint (1994),
to appear in Amer.J.Math.

\entry{Ka-1}{Kato M.:}{Examples on an Extension Problem of Holomorphic Maps
and Holomorphic 1-Dimensional Foliations.} Tokyo\ Journal\ Math. {\bf13}, n 1,
 139-146, (1990).

\entry{Ka-2}{Kato M.:}{Compact quotient manifolds of domains in a complex
3-dimensional projective space and the Lebesgue measure of limit sets.}
 Tokyo\ Journal\ Math. {\bf19}, 99-119 (1996).

 \entry{Ka-3}{Kato M.:}{Compact complex manifolds containing "global"
 spherical shells I.} Proc. Intl. Symp. Algebraic Geometry, Kyoto,
 45-84 (1977).

\entry{Kl}{Klimek M.:}{Pluripotential theory.} London.\ Math.\ Soc.\ Monographs
 ,\ New\ Series 6, (1991).

\entry{Lg}{Lelong P.:}{Plurisubharmonic Functions and Positive Differential
Forms. Gordon and Breach.} New-York. (1969) 78 p.

\entry{Lv}{Levi E.:}{Studii sui punti singolari essenziali delle funzioni
analitiche di due o pi\'u variabili complesse.} Annali di Mat. pura ed appl.
 {\bf17}, n 3, 61-87 (1910).

\entry{Re}{Remmert R.:}{Holomorphe und meromorphe Abbildungen komplexer
R\"aume.} Math.\ Ann. {\bf133}, 328-370, (1957).

\entry{Rs}{Rossi H.:}{Attaching analytic spaces to an analytic space along
a pseudoconvex boundary.} Proc.\ Conf.\ Complex\ Analysis, Monneapolis,
Springer-Verlag, 242-256 (1965).

\bigskip\bigskip
Universit\'e de Lille-I

U.F.R. de Math\'ematiques

Villeneuve d'Ascq Cedex

59655 France

ivachkov@gat.univ-lille1.fr

\bigskip\bigskip
IAPMM Acad Sci. of Ukraine

Naukova 3/b, 290053 Lviv

Ukraine
\end